\documentclass[a4paper,12pt]{article}
\usepackage[utf8]{inputenc}
\usepackage{amsmath}
\usepackage{amssymb}
\usepackage{mathtools}
\usepackage{fancyhdr}
\usepackage{siunitx}
\usepackage{mleftright}
\usepackage{booktabs}
\usepackage{caption}
\usepackage{aligned-overset}
\usepackage{subcaption}
\usepackage[intlimits]{esint}
\usepackage{stackengine}
\usepackage[explicit]{titlesec}
\usepackage{tocloft}
\usepackage{lipsum}
\usepackage{authblk}
\usepackage[textwidth=15cm,textheight=22cm]{geometry}

\renewcommand{\headrulewidth}{0.4pt}

\DeclareMathAlphabet{\mathbbold}{U}{bbold}{m}{n}

\usepackage[hyphens,spaces,obeyspaces]{url}
\usepackage[pdftitle={Operational Methods Applied to the Spherical Mean and X-Ray Transform}, pdfauthor={Julius Lehmann},colorlinks=true, citecolor=blue, urlcolor=black, linkcolor=blue]{hyperref}
\usepackage[style=numeric]{biblatex}
\newcommand{\dr}[1]{\mathrm{#1}}
\newcommand{\dd}[1]{\dr{d}#1\!\mathop{}}
\newcommand{\der}[1]{\frac{\mathrm{d}}{\mathrm{d}#1}}
\newcommand{\Der}[2]{\frac{\mathrm{d}#1}{\mathrm{d}#2}}
\newcommand{\DerN}[3]{\frac{\mathrm{d}^{#3}#1}{\mathrm{d}{#2}^{#3}}}
\newcommand{\pde}[1]{\frac{\partial}{\partial #1}}
\newcommand{\Pde}[2]{\frac{\partial #1}{\partial #2}}

\newcommand{\nexp}[1]{\dr{e}^{#1}}

\newcommand{\x}{\dr{i}}
\newcommand{\Nabla}{\Vec{\nabla}}

\newcommand{\laplace}{\mathop{{}\bigtriangleup}\nolimits\!}

\DeclarePairedDelimiter\klam{(}{)}
\DeclarePairedDelimiter\reklam{[}{]}
\DeclarePairedDelimiter\abs{\lvert}{\rvert}
\DeclarePairedDelimiter\norm{\lVert}{\rVert}

\DeclarePairedDelimiterX\braket[2]{\langle}{\rangle}{#1 \delimsize\vert #2}

\title{Operational Methods Applied to the Spherical Mean and X-Ray Transform}
\author{Julius Lehmann\footnote{e-mail: \texttt{julius.lehmann(at)tum.de}}}
\affil{\textit{Physics of Complex Biosystems, Technical University of Munich, 85748 Garching, Germany}}
\date{February 19, 2024\\[0.2ex]\small(updated: \today)}

\normalsize

\begin{document}

\maketitle

\begin{abstract}
    We employ the framework of operational calculus to derive the operators associated with the spherical mean and a class of related averaging means of a function in $n$-dimensional space. Beginning with the classical definition of the spherical mean, we obtain a compact operator representation in terms of confluent hypergeometric functions of the Laplacian. This operator-based formulation provides a straightforward approach to the analysis of spherical means, allowing us to determine their power series expansions, construct series solutions to the corresponding inversion problems, derive the partial differential equations they satisfy, and give meaning to iterated and fractional spherical means. Finally, we apply the spherical mean operator to derive the inversion formula for the X-ray transform in an operational manner.

    \smallskip

    
\end{abstract}

\section{Spherical Mean and Its Generalizations}
The general solution of the wave equation in odd spatial dimensions is given by Kirchhoff's formula, which is derive by employing spherical means to reduce the higher-dimensional wave equation to the wave equation in $\mathbb{R}$ and then applying d'Alembert's celebrated formula to solve the resulting reduced system~\cite{evans2022partial}. Spherical means also arise naturally in the theory of harmonic functions as one of their defining characteristics; specifically, a function is harmonic if and only if it coincides with its own spherical mean~\cite{netuka1994mean}. The spherical mean of a function $f$ in $n$ dimensions over a sphere of radius $r$ is defined by the integral~\cite{finch2007spherical}
\begin{equation}\label{eq:sphericalMean}
    \Bar{f}_S(r,\Vec{x})=\frac{1}{S_{n-1}}\int_{\mathclap{\norm{\Vec{u}}=1}}\dr{d}^nu\,f(\Vec{x}+r\Vec{u}),
\end{equation}
where $S_{n-1}=2\pi^{\frac{n}{2}}/\Gamma(\frac{n}{2})$ denotes the surface area of the $n$-dimensional unit sphere. Based on this definition, we can already derive numerous properties of the spherical mean, including the differential equation it satisfies, approximations valid for small radii $r$, and the corresponding inverse operation. However, suppose we had an operator that we could apply directly to the function $f(\Vec{x})$ in order to obtain its spherical mean; in that case, the properties of the spherical mean would follow directly from the properties of this operator. In the following, we employ operational calculus~\cite{OperationalCalculus} to identify the operator associated with the spherical mean of a function and to use this operator to derive several properties of Eq.~\eqref{eq:sphericalMean}. In essence, operational calculus provides a framework in which we manipulate operators algebraically and assign precise meaning to functions of operators. Although its origins trace back to the foundations of calculus with Gottfried Wilhelm Leibniz, mathematicians in Britain and Ireland during the 19${}^\text{th}$ century, including Sylvester, Boole, Glaisher, Crofton, and Blizard, substantially developed this approach~\cite{bell1938history,cartier2000mathemagics}, and Oliver Heaviside only fully refined it in the 1890s~\cite{nahin2002oliver}. A well-known example of this approach is the operator we employ the most in this paper, which, as suggested by the structure of the integral in Eq.~\eqref{eq:sphericalMean}, is the shift operator~\cite{jordan1965calculus} and is defined in one dimension by the exponential of the derivative operator:
\begin{equation}
    \nexp{t\pde{x}}f(x)=f(x+t).
\end{equation}
We can readily generalize this expression to higher dimensions by invoking the standard properties of the exponential function together with the commutativity of the individual derivative operators, yielding
\begin{equation}\label{eq:Shift}
    \nexp{\Vec{t}\cdot\Nabla_{\Vec{x}}}f(\Vec{x})=f(\Vec{x}+\Vec{t}\,).
\end{equation}

\begin{figure}
    \centering
    \includegraphics[width=\linewidth]{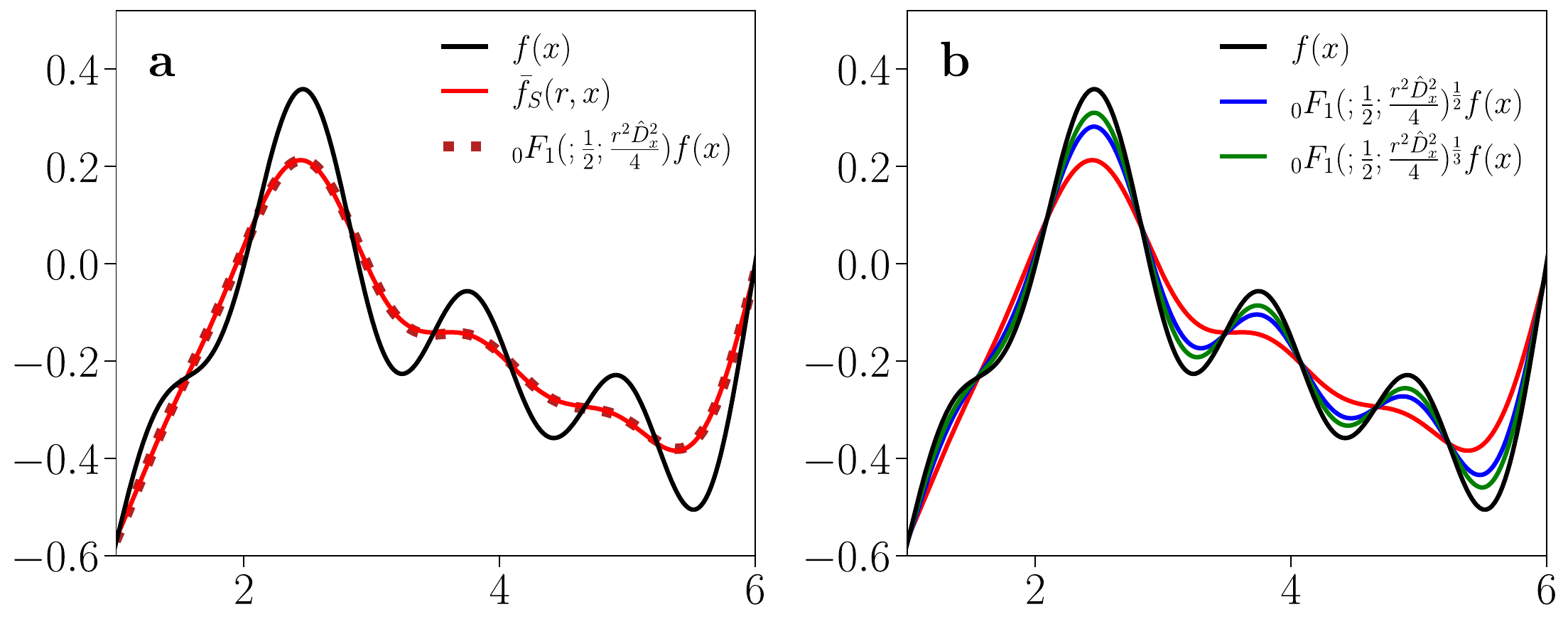}
    \caption{\textbf{Demonstration of different spherical means in 1D.} \textbf{(a)} Comparison of the spherical mean obtained from the integral in Eq.~\eqref{eq:sphericalMean} (solid red line) with the spherical mean obtained from the operator in Eq.~\eqref{eq:sphericalMeanOp} (dashed red line) for the same radius applied to an arbitrary function $f(x)$. \textbf{(b)} Illustration of various fractional spherical means. The solid blue line represents the $\frac{1}{2}$-th spherical mean, and the solid green line represents the $\frac{1}{3}$-rd spherical mean. For reference, the full spherical mean is shown in red.}
    \label{fig:1}
\end{figure}

As mentioned above, we can use this relation to obtain a compact, albeit somewhat unconventional, operator expression for the spherical mean. By directly applying Eq.~\eqref{eq:Shift} to Eq.~\eqref{eq:sphericalMean}, we find
\begin{equation}
    \Bar{f}_S(r,\Vec{x})=\frac{1}{S_{n-1}}\int_{\mathclap{\norm{\Vec{u}}=1}}\dr{d}^nu\, f(\Vec{x}+r\Vec{u})=\frac{1}{S_{n-1}}\int_{\mathclap{\norm{\Vec{u}}=1}}\dr{d}^nu\,\nexp{r\Vec{u}\cdot\Nabla_{\Vec{x}}}f(\Vec{x}).
\end{equation}
The integral on the right-hand side is independent of the specific form of $f(\Vec{x})$ and can be evaluated using standard techniques~\cite{bochner1935summation}, which yields the operator representation
\begin{equation}
    \Bar{f}_S(r,\Vec{x})=\Gamma\klam*{\frac{n}{2}}\klam[\bigg]{\frac{2}{r\norm{\Nabla_{\Vec{x}}}}}^{\frac{n}{2}-1}I_{\frac{n}{2}-1}(r\norm{\Nabla_{\Vec{x}}})f(\Vec{x}),
\end{equation}
where $I_\alpha(x)$ denotes the modified Bessel function of the first kind~\cite{abramowitz1968handbook}. We can express this result in an even more compact form by using generalized hypergeometric functions~\cite{askey2010generalized},
\begin{equation}
    {}_pF_q(a_1,\dots,a_p;b_1,\dots,b_q;x)=\sum_{n=0}^\infty\prod_{i=1}^p\frac{\Gamma(a_i+n)}{\Gamma(a_i)}\prod_{j=1}^q\frac{\Gamma(b_j)}{\Gamma(b_j+n)}\frac{x^n}{n!},
\end{equation}
or, more specifically, the confluent hypergeometric limit function, which allows us to write
\begin{equation}\label{eq:sphericalMeanOp}
    \Bar{f}_S(r,\Vec{x})={}_0F_1\klam[\Big]{;\frac{n}{2};\frac{r^2\laplace_{\Vec{x}}}{4}}f(\Vec{x}).
\end{equation}
This operator is a function of the Laplacian $\laplace_{\Vec{x}}$ acting on $f(\Vec{x})$ and therefore contains all relevant information about the spherical mean of the function. To give this operator a precise meaning, we interpret it through the power series expansion of ${}_0F_1(;\frac{n}{2};x)$ about $x=0$:
\begin{equation}
    {}_0F_1\klam[\Big]{;\frac{n}{2};\frac{r^2\laplace_{\Vec{x}}}{4}}f(\Vec{x})=\klam[\Big]{1+\frac{r^2}{2n}\laplace_{\Vec{x}}+\frac{r^4}{8n(n+2)}\laplace_{\Vec{x}}^2+\dotsi}f(\Vec{x}).
\end{equation}
This expansion immediately provides an efficient means of approximating the spherical mean for small radii $r$. Using standard techniques based on Taylor's theorem, Ref.~\cite{ovall2016laplacian} obtained the same expression, thereby confirming the validity of this result. Before proceeding further, we emphasize that operational methods are not rigorous and must therefore be verified using conventional analytical techniques. Nevertheless, as demonstrated above, they offer valuable insights and practical ``shortcuts'' for deriving analytical results that would otherwise be considerably more involved. With the spherical mean operator at hand, we can readily write down the inverse operation of spherical averaging by dividing both sides by ${}_0F_1\klam{;\frac{n}{2};\frac{r^2\laplace_{\Vec{x}}}{4}}$:
\begin{equation}\label{eq:Inv}
    f(\Vec{x})=\reklam[\Big]{{}_0F_1\klam[\Big]{;\frac{n}{2};\frac{r^2\laplace_{\Vec{x}}}{4}}}^{-1}\Bar{f}_S(r,\Vec{x}).
\end{equation}
Interpreting this inverse operator via its power series expansion, we obtain the corresponding series solution to the inversion problem as
\begin{equation}
    f(\Vec{x})=\klam[\Big]{1-\frac{r^2}{2n}\laplace_{\Vec{x}}+\frac{(n+4)r^4}{8n^2(n+2)}\laplace_{\Vec{x}}^2+\dotsi}\Bar{f}_S(r,\Vec{x}).
\end{equation}
Moreover, we can compute successive spherical means of the same function by applying the operator~\eqref{eq:sphericalMeanOp} repeatedly. Considering the $m$-th iterated spherical mean, we obtain the following power series expansion for small $r$:
\begin{equation}
    \reklam[\Big]{{}_0F_1\klam[\Big]{;\frac{n}{2}; \frac{r^2\laplace_{\Vec{x}}}{4}}}^mf(\Vec{x})=\klam[\Big]{1+\frac{mr^2}{2n}\laplace_{\Vec{x}}+\frac{m(m(n+2)-2)r^4}{8n^2(n+2)}\laplace_{\Vec{x}}^2+\dotsi}f(\Vec{x}).
\end{equation}
This formula remains valid for $m<0$ and thus generalizes Eq.~\eqref{eq:Inv}. By analytic continuation, it also provides a natural framework for defining and computing fractional spherical means, such as the ``$\frac{1}{2}$-th spherical mean'' of $f(\Vec{x})$. In Fig.~\ref{fig:1}, we illustrate several fractional spherical means in one spatial dimension. From the operator formulation, it is also evident that successive spherical means with different radii $r_1$ and $r_2$ commute, that is,
\begin{equation}
    {}_0F_1\klam[\Big]{;\frac{n}{2};\frac{r_1^2\laplace_{\Vec{x}}}{4}}{}_0F_1\klam[\Big]{;\frac{n}{2};\frac{r_2^2\laplace_{\Vec{x}}}{4}}f(\Vec{x})={}_0F_1\klam[\Big]{;\frac{n}{2};\frac{r_2^2\laplace_{\Vec{x}}}{4}}{}_0F_1\klam[\Big]{;\frac{n}{2};\frac{r_1^2\laplace_{\Vec{x}}}{4}}f(\Vec{x}),
\end{equation}
since the underlying operators commute. Nevertheless, spherical means are not additive in the sense that ${}_0F_1\klam[\big]{;\frac{n}{2};\frac{r_2^2\laplace_{\Vec{x}}}{4}}\Bar{f}(r_1,\vec{x})\neq\Bar{f}(r_1+r_2,\vec{x})$. Instead, they satisfy an addition formula involving the Gaussian hypergeometric function ${}_2F_1$:~\cite{erdelyi}
\begin{equation}
    \begin{aligned}
        {}_0F_1\klam[\Big]{;\frac{n}{2};\frac{r_1^2\laplace_{\Vec{x}}}{4}}&{}_0F_1\klam[\Big]{;\frac{n}{2};\frac{r_2^2\laplace_{\Vec{x}}}{4}}f(\Vec{x})\\
        &=\sum_{k=0}^\infty\frac{r_1^{2k}\Gamma(\frac{n}{2})}{4^kk!\Gamma(k+\frac{n}{2})}{}_2F_1\klam[\Big]{1-\frac{n}{2}-k,-k;\frac{n}{2};\klam[\Big]{\frac{r_2}{r_1}}^2}\laplace_{\vec{x}}^kf(\vec{x}).
    \end{aligned}
\end{equation}
In particular, for $r_1=r_2=r$, corresponding to two consecutive spherical means of equal radius, we obtain the comparatively simple hypergeometric representation~\cite{BaileyProductsOG}
\begin{equation}
    \reklam[\Big]{{}_0F_1\klam[\Big]{;\frac{n}{2}; \frac{r^2\laplace_{\Vec{x}}}{4}}}^2f(\Vec{x})={}_1F_2\klam[\Big]{\frac{n-1}{2};\frac{n}{2},n-1; r^2\laplace_{\Vec{x}}}f(\Vec{x}).
\end{equation}

For the next paragraph, we recall another fundamental property of harmonic functions. These functions, which coincide with their own spherical mean, are known to satisfy Laplace's equation,
\begin{equation}
    \laplace_ {\Vec{x}}h(\Vec{x})=0,\quad \Vec{x}\in\Omega,
\end{equation}
subject to appropriate boundary conditions on $\Omega$~\cite{axler2013harmonic}. Although this property does not generally hold for the spherical mean of an arbitrary function, the spherical mean does satisfy a modified version of Laplace's equation that includes additional source terms. To derive the corresponding differential equation, we recall that the confluent hypergeometric limit functions $w\equiv{}_0F_1(;n;x)$ satisfy the ordinary differential equation~\cite{askey2010generalized}
\begin{equation}\label{eq:ODEHyper}
    x\DerN{w}{x}{2}+n\Der{w}{x}-w=0.
\end{equation}
To bring this equation into a form compatible with our operator, we perform the change of variables $x\to ax^2$ and introduce the notation $s\equiv ax^2$, which yields
\begin{equation}
    s\DerN{w(s)}{s}{2}+n\Der{w(s)}{s}-w(s)=0.
\end{equation}
Applying the chain rule to the differential operators, we obtain for the first derivative
\begin{equation}
    \Der{w(ax^2)}{(ax^2)}=\frac{1}{2ax}\Der{w(ax^2)}{x},
\end{equation}
and for the second derivative
\begin{equation}
    \der{(ax^2)}\Der{w(ax^2)}{(ax^2)}
    =\frac{1}{2ax}\klam*{\frac{1}{2ax}\DerN{w(ax^2)}{x}{2}-\frac{1}{2ax^2}\Der{w(ax^2)}{x}}.
\end{equation}
Consequently, under this change of variables, the differential equation~\eqref{eq:ODEHyper} becomes
\begin{equation}
    \DerN{w(ax^2)}{x}{2}+\frac{2n-1}{x}\Der{w(ax^2)}{x}-4aw(ax^2)=0.
\end{equation}
By formally treating the Laplacian as a scalar quantity, we can identify $a\to\laplace_{\Vec{x}}/4$, which implies that the spherical mean operator formally ``satisfies'' the partial differential equation
\begin{equation}
    \Pde{^2}{r^2}{}_0F_1\klam[\Big]{;\frac{n}{2};\frac{r^2\laplace_{\Vec{x}}}{4}}+\frac{n-1}{r}\pde{r}{}_0F_1\klam[\Big]{; \frac{n}{2};\frac{r^2\laplace_{\Vec{x}}}{4}}=\laplace_{\Vec{x}}{}_0F_1\klam[\Big]{;\frac{n}{2};\frac{r^2\laplace_{\Vec{x}}}{4}}.
\end{equation}
Applying this equation from the left to the function $f(\Vec{x})$ and combining the operator with the function to recover the spherical mean, we obtain the well-known Euler--Poisson--Darboux equation~\cite{courant2008methods,zwillinger2021handbook}:
\begin{equation}
    \Pde{^2}{r^2}\Bar{f}_S(r,\Vec{x})+\frac{n-1}{r}\pde{r}\Bar{f}_S(r,\Vec{x})=\laplace_{\Vec{x}}\Bar{f}_S(r,\Vec{x}).
\end{equation}

Having demonstrated the effectiveness and utility of operational calculus and the operator associated with an integral transformation, we now turn to additional examples. In fact, we can regard the spherical mean as a particular instance of a broader class of averaging operations applied to a function over a spherical region of space:
\begin{equation}\label{eq:MeanGen}
    \Bar{f}_K(r,\Vec{x})=\frac{1}{C_K}\int_{\mathclap{\norm{\Vec{u}}\leq1}}\dr{d}^nu\,K(\Vec{u})f(\Vec{x}+r\Vec{u}),
\end{equation}
where $K(\Vec{u})$ denotes the kernel, or weight, of the mean, which we assume to be radially symmetric, and $C_K$,
\begin{equation}
    C_K=\int_{\mathclap{\norm{\Vec{u}}\leq1}}\dr{d}^nu\,K(\Vec{u})=S_{n-1}\int_0^1\dd{u}u^{n-1}K(u),
\end{equation}
is a normalization constant chosen such that averaging the identity function returns the identity. Within this framework, the standard spherical mean corresponds to the choice of kernel $K(\Vec{u})=\delta(\norm{\Vec{u}}-1)$. To determine the operator associated with Eq.~\eqref{eq:MeanGen}, we recall the expression for the Fourier transform of radially symmetric functions~\cite{grafakos2008classical}:
\begin{equation}
    \mathcal{F}\{f(\norm{\Vec{x}})\}(k)=\frac{(2\pi)^{\frac{n}{2}}}{k^{\frac{n}{2}-1}}\int_0^\infty\dd{u}J_{\frac{n}{2}-1}(ku)u^{\frac{n}{2}}f(u).
\end{equation}
Using the shift operator~\eqref{eq:Shift} in the integral representation,
\begin{equation}
    \Bar{f}_K(r,\Vec{x})=\frac{1}{C_K}\int_{\mathclap{\norm{\Vec{u}}\leq1}}\dr{d}^nu\, K(\Vec{u})f(\Vec{x}+r\Vec{u})=\frac{1}{C_K}\int_{\mathclap{\norm{\Vec{u}}\leq1}}\dr{d}^nu\,K(\Vec{u})\nexp{r\Vec{u}\cdot\Nabla_{\Vec{x}}}f(\Vec{x}),
\end{equation}
we observe that the right-hand side is essentially the Fourier transform of $K(\norm{\Vec{u}})H(1-\norm{\Vec{u}})$ after a ``Wick rotation'' in frequency space $k\to\x r$. Substituting the Fourier transform expression into the integral, we obtain
\begin{equation}\label{eq:MeanGenOp}
    \Bar{f}_K(r,\Vec{x})=\frac{(2\pi)^{\frac{n}{2}}}{C_K(r\norm{\Nabla_{\Vec{x}}})^{\frac{n}{2}-1}}\int_0^1\dd{u}I_{\frac{n}{2}-1}(ru\norm{\Nabla_{\Vec{x}}})K(u)u^{\frac{n}{2}}f(\Vec{x}),
\end{equation}
from which we recognize the factor $(r\norm{\Nabla_{\Vec{x}}})^{-\frac{n}{2}+1}I_{\frac{n}{2}-1}(ru\norm{\Nabla_{\Vec{x}}})$ familiar from the ordinary spherical mean. Consequently, we can equivalently express the integral in terms of the confluent hypergeometric limit function:
\begin{equation}
    \Bar{f}_K(r,\Vec{x})=\frac{2\pi^{\frac{n}{2}}}{C_K\Gamma(\frac{n}{2})}\int_0^1\dd{u}K(u)u^{n-1}{}_0F_1\klam[\Big]{;\frac{n}{2};\frac{(ru)^2\laplace_{\Vec{x}}}{4}}f(\Vec{x}),
\end{equation}
which makes it evident that the generalized spherical mean associated with a kernel $K(u)$ is simply a weighted integral of the spherical mean of $f(\Vec{x})$:
\begin{equation}
    \Bar{f}_K(r,\Vec{x})=\frac{2\pi^{\frac{n}{2}}}{C_K\Gamma(\frac{n}{2})}\int_0^1\dd{u}K(u)u^{n-1}\Bar{f}_S(ru,\Vec{x}).
\end{equation}

As a first simple example, we consider the constant kernel $K(\Vec{u})=1$, which yields the average of a function over the entire $n$-ball rather than only over its surface. From Eq.~\eqref{eq:MeanGenOp}, we directly obtain the corresponding operator:
\begin{equation}\label{eq:ballMean}
    \Bar{f}_B(r,\Vec{x})={}_0F_1\klam[\Big]{;\frac{n}{2}+1;\frac{r^2\laplace_{\Vec{x}}}{4}}f(\Vec{x}).
\end{equation}
As before, we can employ this operator to derive small-$r$ approximations or to construct the series solution of the associated inversion problem. Moreover, because its structure closely resembles that of the spherical mean operator, we can immediately write down the partial differential equation satisfied by the $n$-ball mean:
\begin{equation}
    \Pde{^2}{r^2}\Bar{f}_B(r,\Vec{x})+\frac{n+1}{r}\pde{r}\Bar{f}_B(r,\Vec{x})=\laplace_{\Vec{x}}\Bar{f}_B(r,\Vec{x}).
\end{equation}

A more interesting example arises from the general kernel $K(\Vec{u})=(1-\norm{\Vec{u}}^\beta)^\alpha$ with normalization constant $C_K=\frac{\Gamma(\alpha+1)\Gamma(\frac{n}{\beta})}{\beta\Gamma(\frac{n}{\beta}+\alpha+1)}S_{n-1}$. We now examine two special cases. First, we consider the bell-shaped case $\beta = 2$. Substituting this kernel into Eq.~\eqref{eq:MeanGenOp}, we obtain the familiar operator
\begin{equation}
    \Bar{f}_K(r,\Vec{x})={}_0F_1\klam[\Big]{;\frac{n}{2}+\alpha+1;\frac{r^2\laplace_{\Vec{x}}}{4}}f(\Vec{x}),
\end{equation}
which satisfies the partial differential equation
\begin{equation}
    \Pde{^2}{r^2}\Bar{f}_K(r,\Vec{x})+\frac{2\alpha+n+1}{r}\pde{r}\Bar{f}_K(r,\Vec{x})=\laplace_{\Vec{x}}\Bar{f}_K(r,\Vec{x}).
\end{equation}
Both expressions converge to the $n$-ball mean in the limit $\alpha\to0$, as expected. Furthermore, by using the hypergeometric identity ${}_0F_1(;a;x)=\frac{x}{(a+1)a}{}_0F_1(;a+2;x)+{}_0F_1(;a+1;x)$, we obtain the following notable relationship between the spherical average $\Bar{f}_S$, the ball average $\Bar{f}_B$, and the average over a shifted parabola $\Bar{f}_K$ (with $\alpha=1$):
\begin{gather}
    \Bar{f}_S(r,\Vec{x})=\frac{r^2}{n(n+2)}\laplace_{\Vec{x}}\Bar{f}_K(r,\Vec{x})+\Bar{f}_B(r,\Vec{x}).
\end{gather}

The second case, $\beta = 1$, generalizes the triangular weight and leads to the operator
\begin{equation}
    \Bar{f}_K(r,\Vec{x})={}_1F_2\klam[\Big]{\frac{n+1}{2};\frac{n+\alpha+1}{2},\frac{n+\alpha}{2}+1;\frac{r^2\laplace_{\Vec{x}}}{4}}f(\Vec{x}).
\end{equation}
The associated partial differential equation is considerably more involved and given by
\begin{equation}
    \frac{2r}{n+1}\DerN{\Bar{f}_K}{r}{3}+\frac{2(n+\alpha)}{n+1}\DerN{\Bar{f}_K}{r}{2}+\reklam*{\frac{\alpha^2+n(n-1)+(2n-1)\alpha}{(n+1)r}+\frac{r\laplace_{\Vec{x}}}{(n+1)}}\Der{\Bar{f}_K}{r}=\laplace_{\Vec{x}}\Bar{f}_K.
\end{equation}
Once again, these expressions reduce to the $n$-ball mean in the limit $\alpha\to0$.

Before proceeding to the next section, we briefly discuss two additional examples. The first concerns the one-dimensional version of Eq.~\eqref{eq:ballMean}, which coincides with the moving average of a function over a window of width $r$:
\begin{equation}
    \Bar{f}(r,x)=\frac{1}{2r}\int_{\mathclap{-r}}^r\dr{d}u\,f(x+u)=\frac{\sinh(r\hat{D}_x)}{r\hat{D}_x}f(x),
\end{equation}
where, by employing the two integral identities~\cite{bateman_2023_cnd32-h9x80}
\begin{equation}
    \int_0^\infty\dd{u}\cos(\alpha u)\nexp{-su}=\frac{s}{(s^2+\alpha^2)},\qquad\int_0^\infty\dd{u}\abs{\cos(\alpha u)}\nexp{-su}=\frac{s+\alpha\operatorname{csch}(\frac{\pi}{2\alpha}s)}{(s^2+\alpha^2)},
\end{equation}
we obtain the inverse operation $f(x)=\frac{r\hat{D}_x}{\sinh(r\hat{D}_x)}\Bar{f}(r,x)$, which reconstructs the original function from its moving average:
\begin{equation}
    f(x)=\frac{2r^2}{\pi}\int_0^\infty\dd{u}\klam[\Big]{\abs[\Big]{\cos\klam[\Big]{\frac{\pi}{2r} u}}-\cos\klam[\Big]{\frac{\pi}{2r} u}}\klam[\Big]{\Bar{f}'''(r,x-u)+\klam[\Big]{\frac{\pi}{2r}}^2\Bar{f}'(r,x-u)},
\end{equation}
under suitable regularity conditions on $\Bar{f}(r,x)$.

The second example relates to the two-dimensional spherical mean,
\begin{equation}
    \Bar{f}_S(r,\Vec{x})=\frac{1}{2\pi}\int_{\mathclap{\norm{\Vec{u}}=1}}\dr{d}^2u\,f(\Vec{x}+r\Vec{u}),
\end{equation}
by ``complexifying'' the $y$-component of $\Vec{u}$, which leads to the integral
\begin{equation}
    \frac{1}{2\pi}\int_0^{2\pi}\dd{\theta}f(z+r\nexp{\x\theta}).
\end{equation}
Applying the shift operator to the integrand, however, shows that
\begin{equation}
    \frac{1}{2\pi}\int_0^{2\pi}\dd{\theta}\nexp{r\nexp{\x\theta}\hat{D}_z}f(z)=f(z),
\end{equation}
since $\int_0^{2\pi}\dd{\theta}\nexp{r\nexp{\x\theta}}=2\pi$ for $r\neq0$. Consequently, the ``complex'' two-dimensional spherical mean of a function coincides with the function itself, a result known as Gauss's mean value theorem~\cite{weissteinGauss}. Although this argument does not constitute a rigorous proof because it omits the necessary conditions on $f$, it nevertheless provides an intuitive and concise route to this classical conclusion.

\section{X-Ray Transform}
In this section, we apply the operator associated with the spherical mean to derive the inversion formula for the X-ray transform. This integral transform is closely related to the Radon transform and, in two dimensions, coincides with it. In general, we define the X-ray transform of a function $f$ by integrating $f$ along lines in space. Parameterizing a line by $L=\Vec{x}+\ell\vec{\theta}$, where $\Vec{x}$ denotes an initial point and $\vec{\theta}$ specifies the direction of the line, we can write
\begin{equation}
    X\{f\}(\Vec{x},\vec{\theta})=\int_Lf(L)=\int_{\mathclap{-\infty}}^\infty\dr{d}\ell\,f(\Vec{x}+\ell\vec{\theta}).
\end{equation}
To derive the inversion formula, we integrate the left-hand side over the surface of the unit sphere,
\begin{equation}
    \frac{1}{S_{n-1}}\int_{\mathclap{\norm{\vec{\theta}}=1}}\dr{d}^n\theta\,X\{f\}(\Vec{x},\vec{\theta}),
\end{equation}
and then exchange the order of integration between the surface integral and the X-ray transform:
\begin{equation}
    \frac{1}{S_{n-1}}\int_{\mathclap{\norm{\Vec{\theta}}=1}}\dr{d}^n\theta\int_{\mathclap{-\infty}}^\infty\dr{d}\ell\, f(\Vec{x}+\ell\vec{\theta})=\int_{\mathclap{-\infty}}^\infty\dr{d}\ell\,\frac{1}{S_{n-1}}\int_{\mathclap{\norm{\Vec{\theta}}=1}}\dr{d}^n\theta\,f(\Vec{x}+\ell\vec{\theta}).
\end{equation}
The inner integral on the right-hand side is simply the spherical mean of $f(\Vec{x})$. We can therefore rewrite this expression using the operator form~\eqref{eq:sphericalMeanOp} as
\begin{equation}
    \int_{\mathclap{-\infty}}^\infty\dr{d}\ell\,\frac{1}{S_{n-1}}\int_{\mathclap{\norm{\Vec{\theta}}=1}}\dr{d}^n\theta\, f(\Vec{x}+\ell\vec{\theta})=\int_{\mathclap{-\infty}}^\infty\dr{d}\ell\,{}_0F_1\klam[\Big]{; \frac{n}{2};\frac{\ell^2\laplace_{\Vec{x}}}{4}}f(\Vec{x}),
\end{equation}
which evaluates to the operator
\begin{equation}
    \int_{\mathclap{-\infty}}^\infty\dr{d}\ell\,{}_0F_1\klam[\Big]{;\frac{n}{2};\frac{\ell^2\laplace_{\Vec{x}}}{4}}f(\Vec{x})=\frac{\Gamma(\frac{n}{2})}{\Gamma(\frac{n-1}{2})}\sqrt{-\frac{4\pi}{\laplace_{\Vec{x}}}}f(\Vec{x}).
\end{equation}
By rearranging this expression to solve for $f(\Vec{x})$, we obtain the following formal inversion formula:
\begin{equation}
    f(\Vec{x})=\frac{\Gamma(\frac{n-1}{2})}{\sqrt{4\pi}\Gamma(\frac{n}{2})}\sqrt{-\laplace_{\Vec{x}}}\frac{1}{S_{n-1}}\int_{\mathclap{\norm{\vec{\theta}}=1}}\dr{d}^n\theta\,X\{f\}(\Vec{x},\vec{\theta}).
\end{equation}
This formula, however, depends on the interpretation of the operator $\sqrt{-\laplace_{\Vec{x}}}$ and is therefore of limited practical use without further specification. To give meaning to the square root of the Laplacian, we multiply both the numerator and denominator by $\sqrt{-\laplace_{\Vec{x}}}$, which yields
\begin{equation}
    f(\Vec{x})=-\frac{\Gamma(\frac{n-1}{2})}{\sqrt{4\pi}\Gamma(\frac{n}{2})}\frac{\laplace_{\Vec{x}}}{\sqrt{-\laplace_{\Vec{x}}}}\frac{1}{S_{n-1}}\int_{\mathclap{\norm{\vec{\theta}}=1}}\dr{d}^n\theta\, X\{f\}(\Vec{x},\vec{\theta}).
\end{equation}
Here, the resulting operator $1/\sqrt{-\laplace_{\Vec{x}}}$ admits a rigorous definition through its Fourier multiplier, $\mathcal{F}\{1/\sqrt{-\laplace_{\Vec{x}}}f(\Vec{x})\}(\Vec{k})\propto\norm{\Vec{k}}^{-1}\mathcal{F}\{f\}(\Vec{k})$, and is known as the Riesz potential~\cite{riesz1949integrale, landkof1972foundations}. It is given explicitly by the integral representation
\begin{equation}
    \frac{1}{\sqrt{-\laplace_{\Vec{x}}}}f(\Vec{x})=\frac{\Gamma(\frac{n-1}{2})}{2\pi^{\frac{n}{2}}\Gamma(\frac{n}{2})}\int_{\mathbb{R}^n}\dr{d}^nu\,\frac{f(\Vec{u})}{\norm{\Vec{x}-\Vec{u}}^{n-1}}.
\end{equation}
Substituting this expression into the inversion formula, we recover the well-known result~\cite{lai2022single}
\begin{equation}
    f(\Vec{x})=-\frac{\Gamma(\frac{n-1}{2})^2}{8\pi^{n+1}}\int_{\mathbb{R}^n}\dr{d}^nu\int_{\mathclap{\norm{\vec{\theta}}=1}}\dr{d}^n\theta\, \laplace_{\Vec{x}}\klam[\bigg]{\frac{X\{f\}(\Vec{x},\vec{\theta})}{\norm{\Vec{x}-\Vec{u}}^{n-1}}}.
\end{equation}

\clearpage
\lhead{}
\renewcommand{\headrulewidth}{0.0pt}
\printbibliography

\end{document}